\documentclass[number,citesort,seceqn,dvips]{arxbj}
\usepackage{graphicx}

\aid{0}
\volume{17}
\issue{4}
\pubyear{2011}
\firstpage{1344}
\lastpage{1367}
\doi{10.3150/10-BEJ318}

\makeatletter
\newcommand{\nt}{\lfloor nt\rfloor}
\newcommand{\mt}{\lfloor mt\rfloor}
\newcommand{\Dconv}{\stackrel{\mathcal{D}}{\longrightarrow}}

\newcommand{\var}{\operatorname{var}}
\newtheorem{Theorem}{Theorem}[section]
\newtheorem{Lemma}{Lemma}[section]
\newremark{Example}{Example}[section]
\newcommand{\eqref}[1]{(\ref{#1})}
\makeatother

\begin{document}
\begin{frontmatter}

\title{Asymptotics of trimmed CUSUM statistics}
\runtitle{Asymptotics of trimmed CUSUM statistics}

\begin{aug}
\author[1]{\fnms{Istv\'{a}n} \snm{Berkes}\thanksref{1,e1}\ead[label=e1,mark]{berkes@tugraz.at}},
\author[2]{\fnms{Lajos} \snm{Horv\'{a}th}\thanksref{2}\ead[label=e2]{horvath@math.utah.edu}}
\and
\author[1]{\fnms{Johannes} \snm{Schauer}\thanksref{1,e3}\ead[label=e3,mark]{johannes@schauer.com}\corref{}}

\runauthor{I. Berkes, L. Horv\'{a}th and J. Schauer}
\address[1]{Institute of Statistics, Graz University of
Technology, M\"{u}nzgrabenstra{\ss}e 11, A-8010 Graz, Austria. \mbox{\printead{e1,e3}}}
\address[2]{Department of Mathematics, University of Utah, 155 South
1440 East, Salt Lake City, UT 84112-0090 USA. \printead{e2}}
\end{aug}

\received{\smonth{11} \syear{2009}}
\revised{\smonth{6} \syear{2010}}

\begin{abstract}
There is a wide literature on change point tests, but the case of
variables with infinite variances is essentially unexplored. In this
paper we address this problem by studying the asymptotic behavior of
trimmed CUSUM statistics. We show that in a location model with
i.i.d. errors in the domain of attraction of a stable law of
parameter $0 < \alpha < 2$, the appropriately trimmed CUSUM process
converges weakly to a Brownian bridge. Thus, after moderate trimming,
the classical method for detecting change points remains valid also
for populations with infinite variance. We note that according to the
classical theory, the partial sums of trimmed variables are
generally not asymptotically normal and using random centering in
the test statistics is crucial in the infinite variance case. We
also show that the partial sums of truncated and trimmed random
variables have different asymptotic behavior. Finally, we discuss
resampling procedures which enable one to determine critical values in
the case of small and moderate sample sizes.
\end{abstract}


\end{frontmatter}

\section{Introduction}\label{se:intro}

In this paper we are interested in detecting a possible change in the
location of independent observations. We observe $X_1,\ldots,X_n$ and
want to test the no change null hypothesis
\[
\begin{tabular}{@{}p{331pt}@{}}
$H_0\dvtx X_1,X_2,\ldots,X_n$ are independent, identically
distributed random variables
\end{tabular}
\]
against the $r$ changes alternative\vspace*{-1pt}
\[
H_A\dvtx  X_j= \cases{
e_j, &\quad$1\le j \le n_1$,\cr
e_j+c_1, &\quad$n_1<j\le n_2$,\cr
e_j+c_2, &\quad$n_2<j\le n_3$,\cr
\vdots &\quad\cr
e_j+c_r, &\quad$n_r<j\le n$.}
\]
It is assumed that
\begin{equation}\label{eq:iid}
\begin{tabular}{@{}p{295pt}@{}}
$e_1,\ldots,e_n$ are independent, identically distributed random
variables,
\end{tabular}
\end{equation}
that $c_0=0$, $c_i\neq c_{i+1}$, $i=0,\ldots,r-1$, and that
$1\le n_1<n_2<\cdots<n_r<n$ are unknown. In our model, the changes are at
time $n_j$, $1\le j \le r$. Testing $H_0$ against $H_A$ has been
considered by several authors. For surveys, we refer to
Brodsky and Darkhovsky \cite{broddark1}, Chen and Gupta \cite{chengupta} and Cs\"{o}rg\H{o} and H\'{o}rvath \cite{ch}. If the observations have
finite expected value, then the model is referred to as \textit{changes
in the mean}.\vadjust{\goodbreak}

Several of the most popular methods are based on the functionals of the
CUSUM process (tied down partial sums)
\[
M_n(t)=\sum_{j=1}^{\nt} X_j -\frac{\nt}{n} \sum_{j=1}^n
X_j.
\]
If $H_0$ holds and $0<\sigma^2=\var X_1<\infty$, then
\begin{equation}\label{eq:CUSUM_conv}
\frac{1}{\sqrt{n}} M_n(t)\stackrel{\mathcal{D}[0,1]}{\longrightarrow}
\sigma B(t),
\end{equation}
where $\{B(t),0\le t\le 1\}$ is a Brownian bridge. If $\hat{\sigma}_n$
is a weakly consistent estimator for $\sigma$, that is,
$\hat{\sigma}_n\to \sigma$ in probability, then
\begin{equation}\label{eq:CUSUM_conv2}
\frac{1}{\hat{\sigma}_n \sqrt{n}} M_n(t)
\stackrel{\mathcal{D}[0,1]}{\longrightarrow} B(t).
\end{equation}
Functionals of \eqref{eq:CUSUM_conv2} can be used to find
asymptotically distribution-free procedures to test $H_0$ against
$H_A$. The limit results in \eqref{eq:CUSUM_conv} and
\eqref{eq:CUSUM_conv2} have been extended in several directions. Due to
applications in economics, finance, meteorology, environmental sciences
and quality control, several authors have studied the properties of
$M_n(t)$ and especially \eqref{eq:CUSUM_conv2} for dependent
observations. For relevant references, we refer to
Horv\'{a}th and Steinebach \cite{horvathsteinebach}. The case of vector-valued dependent
observations is considered in Horv\'{a}th, Kokoszka and Steinebach \cite{horkoksteine}. We note that in the
case of dependent observations, $\sigma^2=\lim_{n\to\infty}
\var(n^{-1/2} \sum_{j=1}^n X_j)$, so the estimation of~$\sigma$ is
considerably harder than in the i.i.d. case (see
Bartlett \cite{bartlett}, Grenander and Rosenblatt \cite{grenrosen} and Parzen \cite{parzen}). The rate of convergence in
\eqref{eq:CUSUM_conv2} may be slow, so the asymptotic critical values
might be misleading; hence, resampling methods have been advocated in
Hu\v{s}kov\'{a} \cite{huskova}. With very few exceptions, it has been assumed that at
least $E X_j^2$ is finite. In this paper we are interested in testing
$H_0$ against $H_A$ when $E X_j^2=\infty$.

We assume that
\begin{equation}\label{eq:ass_stable}
\begin{tabular}{@{}p{336pt}@{}}
$X_1,X_2,\ldots$  belong to the domain of attraction of a stable random variable
$\xi_{\alpha}$\\
with parameter $0<\alpha<2$
\end{tabular}
\end{equation}
and
\begin{equation}\label{eq:ass_symm}
X_j \mbox{ is symmetric when }\alpha=1.
\end{equation}
This means that
\begin{equation}\label{eq:stable}
\Biggl(\sum_{j=1}^n X_j-a_n\Biggr)\Bigl/b_n\Dconv \xi_{\alpha}
\end{equation}
for some numerical sequences $a_n$ and $b_n$. The necessary and
sufficient conditions for \eqref{eq:stable} are
\begin{equation}\label{eq:tails_stable}
\lim_{t\to\infty} \frac{P\{X_1>t\}}{L(t) t^{-\alpha}}=p
\quad\mbox{and}\quad
 \lim_{t\to\infty} \frac{P\{X_1\le -t\}}{L(t)t^{-\alpha}}=q
\end{equation}
for some numbers $p\ge 0$, $q\ge 0$ with $p+q=1$ and where $L$ is a
slowly varying function at $\infty$.\vadjust{\goodbreak}

Aue \textit{et al.} \cite{aubeho} studied the properties of $M_n(t)$
under conditions $H_0$, \eqref{eq:ass_stable} and \eqref{eq:ass_symm}.
They used $\max_{1\le j\le n}|X_j|$ as the normalization of $M_n(t)$
and showed that
\begin{equation}\label{eq:CUSUM_conv_max}
\frac{1}{\gamma_n} M_n(t) \stackrel{\mathcal{D}[0,1]}{\longrightarrow}
\frac{1}{\mathcal{Z}} B_{\alpha}(t), \qquad \gamma_n=\max_{1\le j\le
n}|X_j|.
\end{equation}
Here, $B_{\alpha}(t)=W_{\alpha}(t)-t W_{\alpha}(1)$ is an
$\alpha$-stable bridge, $W_{\alpha}(t)$ is an $\alpha$-stable process
(see also Kasahara and Watanabe \cite{kasaharawatanabe}, Section 9) and  $\mathcal{Z}$ is a
random norming factor whose joint distribution with $W_{\alpha}(t)$ is
described in \cite{aubeho} explicitly. Nothing is known about the
distribution of the functionals of $B_{\alpha}(t)/\mathcal{Z}$ and
therefore it is nearly impossible to determine critical values needed
to construct asymptotic test procedures. Hence, resampling methods
(bootstrap and permutation) have been tried. However, it was proven
that the conditional distribution of the resampled $M_n(t)/\gamma_n$,
given $X_1,\ldots,X_n$, converges in distribution to a non-degenerate
random process depending also on the trajectory $(X_1, X_2, \ldots)$.
So, resampling cannot be recommended to obtain asymptotic critical
values. This result was obtained by Aue \textit{et al}. \cite{aubeho}
for permutation resampling and by Athreya \cite{at}, Hall \cite{ha} and
Berkes \textit{et al}. \cite{berhorsch} for the bootstrap. No efficient
procedure has been found to test $H_0$ against $H_A$ when $E
X_j^2=\infty$.

The reason for the `bad' behavior of the CUSUM statistics described
above is the influence of the large elements of the sample. It is known
that for i.i.d. random variables $X_1, X_2, \ldots$ in the domain of
attraction of a non-normal stable law, the $j$th largest element of
$|X_1|, \ldots, |X_n|$ has, for any fixed $j$, the same order of
magnitude as the sum $S_n=X_1+\cdots +X_n$ as $n\to\infty$. Thus, the
influence of the large elements in the CUSUM functional does not become
negligible as $n\to\infty$ and, consequently, the limiting behavior of
the CUSUM statistics along different trajectories $(X_1, X_2, \ldots)$
is different, rendering this statistics impractical  for statistical
inference. The natural remedy for this trouble is trimming, that is,
removing the $d(n)$ elements  with the largest absolute values from the
sample, where $d(n)$ is a suitable number with $d(n)\to\infty$,
$d(n)/n\to 0.$ This type of trimming is usually called \textit{modulus
trimming} in the literature. In another type of trimming, some of the
largest and smallest order statistics are removed from the sample (see,
e.g., Cs\"{o}rg\H{o} et al. \cite{cshaema2,cshormas}). Under suitable conditions, trimming
indeed leads  to a better asymptotic behavior of partial sums (see,
e.g., Mori \cite{mori76,mori77,mori84}, Maller  \cite{maller82,maller84}, Cs\"{o}rg\H{o} et al. \cite{cshaema1,cshaema2,cshormas},
Griffin and Pruitt \cite{griffinpruitt87,griffinpruitt89} and Haeusler and Mason \cite{hama,hama2}). Note, however, that
the asymptotic properties of trimmed random variables depend strongly
on the type of trimming used. In this paper, trimming means modulus
trimming, as introduced above. Griffin and Pruitt
\cite{griffinpruitt87} showed that even in the case where the $X_j$
belong to the domain of attraction of a symmetric stable law with
parameter $0<\alpha<2$, the modulus trimmed partial sums need not be
asymptotically normal. Theorem \ref{th:diff_trim_trunc} reveals the
reason for this surprising fact: for non-symmetric distributions $F$,
the center of the sample remains, even after modulus trimming, a
non-degenerate random variable, and no non-random centering can lead to
a central limit theorem. In contrast, a suitable random centering will
always work and since the CUSUM functional is not affected by centering
factors, even in the case of `bad' partial sum behavior, the trimmed
CUSUM functional converges to a Brownian bridge, resulting in a simple
and useful change point test.

To formulate our results, consider the trimmed CUSUM process
\[
T_n(t)=\sum_{j=1}^{\nt} X_j I\{|X_j|\le \eta_{n,d}\}-\frac{\nt}{n}
\sum_{j=1}^n X_j I\{|X_j|\le \eta_{n,d}\},\qquad
 0\le t\le 1,
 \]
where $\eta_{n,d}$ is the $d$th largest value among
$|X_1|,\ldots,|X_n|$.

Let
\[
F(t)=P\{X_1\le t\}  \quad\mbox{and}\quad  H(t)=P\{|X_1|>t\}.
\]
The (generalized) inverse (or quantile) of $H$ is denoted  $H^{-1}(t)$. We
assume that
\begin{equation}\label{eq:dn}
\lim_{n\to\infty} d(n)/n = 0
\end{equation}
and
\begin{equation}\label{eq:dn2}
\lim_{n\to\infty} d(n)/(\log n)^{7+\varepsilon}=\infty
\qquad\mbox{with some }\varepsilon>0.
\end{equation}
For the sake of simplicity (see Mori \cite{mori76}), we also require that
\begin{equation}\label{eq:F_cont}
F \mbox{ is continuous.}
\end{equation}
Let
\begin{equation}\label{pityu}
A_n^2=\frac{\alpha}{2-\alpha} \bigl(H^{-1}(d/n)\bigr)^2 d.
\end{equation}
Our first result states the weak convergence of $T_n(t)/A_n$.
\begin{Theorem}\label{th:weak_conv_T}
If $H_0$, \eqref{eq:ass_stable}, \eqref{eq:ass_symm} and
\eqref{eq:dn}--\eqref{eq:F_cont} hold, then
\begin{equation}\label{eq:weak_conv_T}
\frac{1}{A_n} T_n(t)\stackrel{\mathcal{D}[0,1]}{\longrightarrow} B(t),
\end{equation}
where $\{B(t),0\le t\le 1\}$ is a Brownian bridge.
\end{Theorem}

Since $A_n$ is unknown, we need to estimate it from the sample. We will
use
\[
\hat{A}_n^2=\sum_{j=1}^n (X_j I\{|X_j|\le
\eta_{n,d}\}-\bar{X}_{n,d})^2\quad\mbox{and}\quad \hat{\sigma}_n^2=
\frac{1}{n}\hat{A}_n^2,
\]
where
\[
\bar{X}_{n,d}=\frac{1}{n} \sum_{j=1}^n X_j I\{|X_j|\le \eta_{n,d}\}.
\]
We note that $\hat{A}_n/A_n\to 1$ almost surely (see
Lemma~\ref{le:2.7}).

\begin{Theorem}\label{th:weak_conv_T2}
If the conditions of Theorem~\ref{th:weak_conv_T} are satisfied, then
\begin{equation}\label{eq:weak_conv_T2}
\frac{1}{\hat{\sigma}_n\sqrt{n}}T_n(t)\stackrel{\mathcal{D}[0,1]}{\longrightarrow}
B(t).
\end{equation}
\end{Theorem}

In the case of independence and $0<\sigma^2=\var X_j<\infty$, we
estimate $\sigma^2$ by the sample variance. So, the comparison of
\eqref{eq:CUSUM_conv2} and \eqref{eq:weak_conv_T2} reveals that in case
of $E X_j^2=\infty,$ we still use the classical CUSUM procedure; only
the extremes are removed from the sample. The finite-sample properties
of tests for $H_0$ against $H_A$ based on \eqref{eq:weak_conv_T2} are
investigated in Section~\ref{se:simulations}.

In the case of a given sample, it is difficult to decide if $E X_j^2$
is finite or infinite. Thus, for applications, it is important to
establish Theorem~\ref{th:weak_conv_T2} when $E X_j^2<\infty$.

\begin{Theorem}\label{th:weak_conv_T3}
If $H_0$, \eqref{eq:dn}, \eqref{eq:dn2} and $E X_j^2<\infty$ are
satisfied, then \eqref{eq:weak_conv_T2} holds.
\end{Theorem}

Combining Theorems~\ref{th:weak_conv_T2} and \ref{th:weak_conv_T3}, we
see that the CUSUM-based procedures can always be used if the
observations with the largest absolute values are removed from the
sample.

We now outline the basic idea of the proofs of
Theorems~\ref{th:weak_conv_T} and \ref{th:weak_conv_T2}. It was proven
by Kiefer~\cite{kiefer71} (see Shorack and Wellner \cite{showe}) that $\eta_{n,d}$ is close
to $H^{-1}(d/n)$ and thus it is natural to consider the process
obtained from $T_n(t)$ by replacing $\eta_{n, d}$ with $H^{-1}(d/n)$.
Let
\[
V_n(t)=\sum_{j=1}^{\nt} \bigl(X_j I\{|X_j|\le H^{-1}(d/n)\} -E\bigl(X_j
I\{|X_j|\le H^{-1}(d/n)\}\bigr)\bigr)
\]
and
\[
V_n^*(t)=\sum_{j=1}^{\nt} \bigl(X_j I\{|X_j|\le \eta_{n, d}\} -E(X_j
I\{|X_j|\le \eta_{n, d}\})\bigr).
\]
Since $V_n(t)$ is a sum of i.i.d. random variables, the classical
functional central limit theorem for triangular arrays easily yields
the following result.

\begin{Theorem}\label{th:weak_conv_V} If the
conditions of Theorem~\ref{th:weak_conv_T} are satisfied, then
\[
\frac{1}{A_n} V_n(t)\stackrel{\mathcal{D}[0,1]}{\longrightarrow} W(t),
\]
where $\{W(t),0\le t\le 1\}$ is a standard Brownian motion (Wiener
process).
\end{Theorem}

In view of the closeness of $\eta_{n,d}$ and $H^{-1}(d/n)$,  one would
 expect  the asymptotic behavior of $V_n(t)/A_n$ and $V_n^*(t)/A_n$
to be the same. Surprisingly, this is not the case. Let
\[
m(t)=E[X_1 I\{|X_1|\le t\}-X_1 I\{|X_1|\le H^{-1} (d/n)\}], \qquad t\ge
0.
\]

\begin{Theorem}\label{th:diff_trim_trunc}
If the conditions of Theorem~\ref{th:weak_conv_T} are satisfied, then
\[
\frac{1}{A_n} \max_{1\le k\le n} \Biggl|\sum_{j=1}^k \bigl[X_j
\bigl(I\{|X_j|\le \eta_{n,d}\}-I\{|X_j|\le
H^{-1}(d/n)\}\bigr)-m(\eta_{n,d})\bigr] \Biggr|=\mathrm{o}_P(1).
\]
\end{Theorem}

By Theorem~\ref{th:diff_trim_trunc}, the asymptotic properties of the
partial sums of trimmed and truncated variables would be the same if
$n|m(\eta_{n,d})|=\mathrm{o}_P(A_n)$ were true. However, this is not
always the case, as the following example shows.

\begin{Example}\label{ex:diff_trim_trunc}
Assume that $X_1$ is concentrated on $(0, +\infty)$ and has a
continuous  density $f$ which is regularly varying at $\infty$ with
exponent $-(\alpha+1)$ for some $0<\alpha<2.$  Then,
\[
\frac{nm(\eta_{n,d})}{B_n}\Dconv N(0, 1),
\]
where
\[
B_n=\frac{\alpha d^{3/2}}{nH'(H^{-1}(d/n))}.
\]
\end{Example}

We conjecture that the centering factor $n m(\eta_{n,d})/A_n$ and the
partial sum process
\[
\sum_{j=1}^{\nt} \bigl(X_j I\{|X_j|\le H^{-1}(d/n)\}-E\bigl(X_j
I\{|X_j|\le H^{-1}(d/n)\}\bigr)\bigr), \qquad 0\le t\le 1,
\]
are asymptotically independent under the conditions of
Example~\ref{ex:diff_trim_trunc}. Hence, by
Theorem~\ref{th:diff_trim_trunc} we would have
\[
\frac{1}{A_n}\sum_{j=1}^{\nt}(X_j I\{|X_j|\le
\eta_{n,d}\}-c_n)\stackrel{\mathcal{D}[0,1]}{\longrightarrow}W(t)+
t\biggl(\frac{2-\alpha}{\alpha}\biggr)^{1/2} \xi,
\]
where $\{W(t),0\le t\le 1\}$ and $\xi$ are independent, $W(t)$ is a
standard Wiener process, $\xi$ is a standard normal random variable and
$c_n=E X_1 I\{|X_1|\le H^{-1}(d/n)\}$.

In view of Theorem \ref{th:diff_trim_trunc}, the normed partial sum
processes of  $X_j I\{|X_j|\le \eta_{n,d}\}-m(\eta_{n, d})$ and $X_j
I\{|X_j|\le H^{-1}(d/n)\}$ have the same asymptotic behavior and thus
the same holds for the corresponding CUSUM processes. By Theorem
\ref{th:weak_conv_V}, the CUSUM process of $X_j I\{|X_j|\le
H^{-1}(d/n)\}$ converges weakly to the Brownian bridge and the CUSUM
process of $X_j I\{|X_j|\le \eta_{n,d}\}-m(\eta_{n, d})$ clearly
remains the same if we drop the term $m(\eta_{n, d})$. Formally,
\begin{eqnarray}\label{eq:star}
&&\max_{1\le k\le n}\Biggl|\sum_{j=1}^k X_j I\{|X_j|\le
\eta_{n,d}\}-\frac{k}{n} \sum_{j=1}^n X_j I\{|X_j|\le
\eta_{n,d}\} \nonumber\\
&&\quad\qquad-\Biggl(\sum_{j=1}^k X_j I\{|X_j|\le
H^{-1}(d/n)\}-\frac{k}{n}\sum_{j=1}^n X_j I\{|X_j|\le
H^{-1}(d/n)\}\Biggr)\Biggr|\\
&&\quad\le 2\max_{1\le k\le n} \Biggl|\sum_{j=1}^k \bigl[X_j
\bigl(I\{|X_j|\le \eta_{n,d}\}-I\{|X_j|\le
H^{-1}(d/n)\}\bigr)-m(\eta_{n,d})\bigr]\Biggr|.\nonumber
\end{eqnarray}
Thus, even though the partial sums of trimmed and truncated variables
are asymptotically different due to the presence of the random
centering $m(\eta_{n, d})$, the asymptotic distributions of the CUSUM
processes of the trimmed and truncated variables are the same.

The proofs of the asymptotic results for  $\sum_{j=1}^n X_j I\{|X_j|\le
\eta_{n,d}\}$ in
Griffin and Pruitt \cite{griffinpruitt87,griffinpruitt89}, Maller \cite{maller82,maller84}, Mori \cite{mori76,mori77,mori84}
are based on classical probability theory. Cs\"{o}rg\H{o} \textit{et
al.} \cite{cshaema1,cshaema2,cshormas} and Haeusler and Mason
\cite{hama} use the weighted approximation of quantile processes
to establish the normality of a class of trimmed partial sums. The
method of our paper is completely different. We show in
Theorem~\ref{th:diff_trim_trunc} that after a suitable random
centering, trimmed partial sums can be replaced with truncated ones,
reducing the problem to sums of i.i.d.~random variables.

\section{Resampling methods}\label{se:resampling}

Since the convergence in Theorem \ref{th:weak_conv_T} can be slow,
critical values in the change point test determined on the basis of the
limit distribution may not be appropriate for small sample sizes. To
resolve this difficulty, resampling methods can be used to simulate
critical values. Let
\[
x_j=X_j I\{|X_j|\le \eta_{n,d}\}-\bar{X}_{n,d},\qquad 1\le j\le n,
\]
be the trimmed and centered observations. We select $m$ elements from
the set $\{x_1,x_2,\ldots,x_n\}$ randomly (with or without
replacement), resulting in the sample $y_1,\ldots,y_m$. If we select
with replacement, the procedure is the bootstrap; if we select without
replacement and $m=n$, this is the permutation method (see
Hu\v{s}kov\'{a} \cite{huskova}). We now define the resampled CUSUM process
\[
T_{m,n}(t)=\sum_{j=1}^{\mt} y_j
 -\frac{\mt}{m} \sum_{j=1}^{m} y_j.
\]
We note that, conditionally on $X_1,X_2,\ldots ,X_n$, the mean of $y_j$
is 0 and its variance is $\hat{\sigma}^2_n$.

\begin{Theorem}\label{th:resampling}
Assume that the conditions of Theorem~\ref{th:weak_conv_T} are
satisfied and draw $m=m(n)$ elements $y_1,\ldots,y_m$ from the set
$\{x_1,\ldots,x_n\}$ with or without replacement, where
\begin{equation}\label{eq:mn}
m=m(n)\rightarrow \infty\qquad\mbox{as }n\rightarrow \infty
\end{equation}
and $m(n)\le n$ in case of selection without replacement. Then, for
almost all realizations of $X_1,X_2,\ldots,$ we have
\[
\frac{1}{\hat{\sigma}_n\sqrt{m}}
T_{m,n}(t)\stackrel{\mathcal{D}[0,1]}{\longrightarrow}B(t),
\]
where $\{B(t),0\le t\le 1\}$ is a Brownian bridge.
\end{Theorem}

By the results of Aue \textit{et al.} \cite{aubeho} and Berkes
\textit{et al.} \cite{berhorsch}, if we sample from the original
(untrimmed) observations, then the CUSUM process converges weakly to a
non-Gaussian process containing random parameters and thus the
resampling procedure is statistically useless.

If we use resampling to determine critical values in the CUSUM test, we
need to study the limit also under the the alternative since in a
practical situation we do not know which of $H_0$ or $H_A$ is valid. As
before, we assume that the error terms $\{e_j\}$ are in the domain of
attraction of a stable law, that is,
\begin{equation}\label{eq:e_i_stable}
\lim_{t\to\infty} \frac{P\{e_1>t\}}{L(t)t^{-\alpha}}=p \quad
\mbox{and}\quad
 \lim_{t\to\infty} \frac{P\{e_1\le -t\}}{L(t)t^{-\alpha}}=q,
\end{equation}
where $p\ge 0$, $q\ge 0$, $p+q=1$ and $L$ is a slowly varying function
at $\infty$.

\begin{Theorem}\label{th:resampling_H_A}
If $H_A$, \eqref{eq:iid}, \eqref{eq:dn}--\eqref{eq:F_cont},
\eqref{eq:mn} and \eqref{eq:e_i_stable} hold, then for almost all
realizations of $X_1,X_2,\ldots,$ we have that
\[
\frac{1}{\hat{\sigma}_n\sqrt{m}}
T_{m,n}(t)\stackrel{\mathcal{D}[0,1]}{\longrightarrow} B(t),
\]
where $\{B(t),0\le t\le 1\}$ is a Brownian bridge.
\end{Theorem}

In other words, the limiting distribution of the trimmed CUSUM process
is the same under $H_0$ and $H_A$, and thus the critical values
determined by resampling will always work. On the other hand, under
$H_A$, the test statistic $\sup_{0<t<1} |T_n(t)|/A_n$ goes to infinity,
so using the critical values determined by resampling, we get a
consistent test.

\begin{table}[b]
\tablewidth=230pt
\caption{Simulated critical values of $\sup_{0<t<1}
|T_{n}(t)|/(\hat{\sigma}_n \sqrt{n})$ for
$1-\alpha=0.95$}\label{ta:sim_crit}
\begin{tabular*}{230pt}{@{\extracolsep{\fill}}lllll@{}}
\hline
$n=100$ & $n=200$ & $n=400$ & $n=800$ & $n=\infty$\\
\hline
1.244 & 1.272 & 1.299 & 1.312 & 1.358\\
\hline
\end{tabular*}
\end{table}

\begin{figure}[b]

\includegraphics{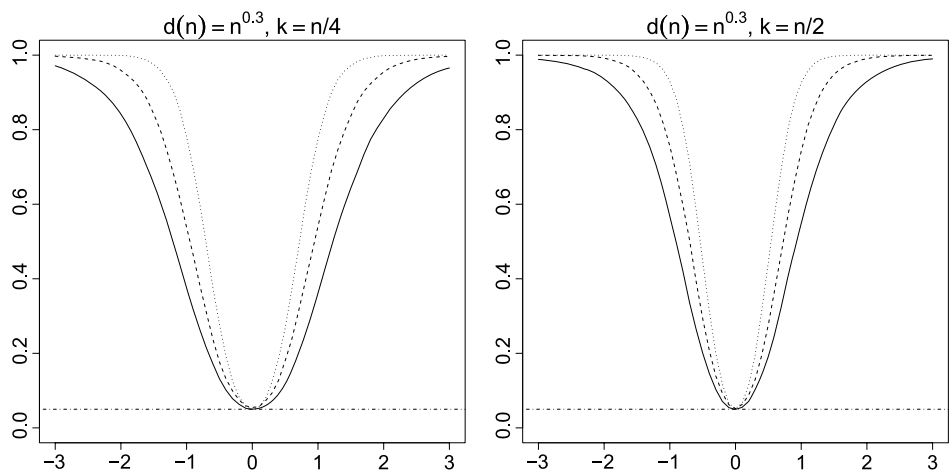}

\caption{Empirical power curves with $\alpha=0.05$,
$n=100$ (solid), $n=200$ (dashed) and $n=400$ (dotted).}
\label{fig:trim_power_1}
\end{figure}

We note that Theorems~\ref{th:resampling} and \ref{th:resampling_H_A}
remain true if \eqref{eq:stable} is replaced with $E X_j^2<\infty$. The
proofs are similar to that of Theorem~\ref{th:resampling} but much
simpler, so no details are given.

\section{Simulation study}\label{se:simulations}

Consider the model under $H_0$ with i.i.d. random variables $X_j$,
$j=1,\ldots,n$, having distribution function
\[
F(t)=\cases{
q(1-t)^{-1.5} &\quad for $t\le 0,$\cr
1-p (1+t)^{-1.5} &\quad for $t > 0,$}
\]
where $p\ge 0$, $q\ge 0$ and $p+q=1$. We trim the samples using
$d(n)=\lfloor n^{0.3} \rfloor$. To simulate the critical values, we
generate $N=10^5$ Monte Carlo simulations for each
$n\in\{100,200,400,800\}$ according to the model under the no change
hypothesis and calculate the values of $\sup_{0<t<1}
|T_{n}(t)|/(\hat{\sigma}_n\sqrt{n})$, where $T_{n}(t)$ and
$\hat{\sigma}_n$ are defined in Section~\ref{se:intro}. The computation
of the empirical quantiles yields the estimated critical values. Table
\ref{ta:sim_crit} summarizes the results for $p=q=1/2$ and
$1-\alpha=0.95$.

Figure~\ref{fig:trim_power_1} shows the empirical power of the test of
$H_0$ against $H_A$ based on the statistic $\sup_{0<t<1}
|T_n(t)|/(\hat{\sigma}_n \sqrt{n})$ for a single change at time
$k=n_1\in\{n/4,n/2\}$ and each $c_1 \in \{-3,-2.9,\ldots,2.9,3\}$
for the same trimming as above ($d(n)=\lfloor n^{0.3} \rfloor$) and a
significance level of $1-\alpha=0.95$, where the number of repetitions
is $N=10^4$. Note that depending on the sample size, we used different
simulated quantiles (see Table~\ref{ta:sim_crit}). The power behaves
best for a change point in the middle of the observation period
($k=n/2$). Due to the differences between the simulated and asymptotic
critical values in Table~\ref{ta:sim_crit}, especially for small $n$,
the test based on the asymptotic critical values tends to be
conservative.

\section{Proofs}\label{se:proofs}

Throughout this section we assume that $H_0$ holds. Clearly,
\[
H(x)=1-F(x)+F(-x),\qquad x\ge 0,
\]
and by \eqref{eq:tails_stable}, we have that
\begin{equation}\label{eq:H_inv}
H^{-1}(t)=t^{-1/\alpha} K(t), \qquad \mbox{if }t\le t_0,
\end{equation}
where $K(t)$ is a slowly varying function at 0. We also use
\begin{equation}\label{eq:d_infty}
d=d(n)\to\infty.
\end{equation}

\begin{Lemma}\label{le:2.1}
If $H_0$, \eqref{eq:ass_stable}, \eqref{eq:ass_symm}, \eqref{eq:dn} and
\eqref{eq:d_infty} hold, then
\begin{equation}\label{eq:V_var}
\lim_{n\to\infty} \frac{1}{A_n^2} \var V_n (1)=1
\end{equation}
and
\begin{eqnarray}\label{eq:4_mom}
&&\lim_{n\to\infty} \sum_{j=1}^n E\bigl[ X_j I\{|X_j|\le
H^{-1}(d/n)\}-E[X_j I\{|X_j|\le H^{-1}(d/n)\}]\bigl]^4\nonumber\\[-8pt]\\[-8pt]
&&\quad{}\times \frac{1}{d(H^{-1}(d/n))^4}
 =\frac{\alpha}{4-\alpha}.\nonumber
\end{eqnarray}
\end{Lemma}

\begin{pf}
If $1<\alpha<2$, then
\[
\lim_{n\to\infty} EX_1 I\{|X_1|\le H^{-1}(d/n)\}=EX_1.
\]
If $\alpha=1$, then by the assumed symmetry, $EX_1 I\{|X_1|\le
H^{-1}(d/n)\}=0$. In the case $0<\alpha<1,$ we write
\begin{eqnarray*}
E|X_1| I\{|X_1|\le H^{-1}(d/n)\} &=&\int_{-H^{-1}(d/n)}^{H^{-1}(d/n)}
|x|\,\mathrm{d}F(x)\\
&=&-\int_{0}^{H^{-1}(d/n)} x \,\mathrm{d}H(x)\\
&=&-x H(x)|_{H^{-1}(d/n)} + \int_{0}^{H^{-1}(d/n)} H(x) \,\mathrm{d}x.
\end{eqnarray*}
By Bingham \textit{et al}. \cite{binghametal}, page~26,
\[
\lim_{y\to\infty} \frac{\int_{0}^{y} H(x)\,
\mathrm{d}x}{(1/(1-\alpha))y^{1-\alpha} L(y)}=1
\]
and therefore
\[
\lim_{n\to\infty} \frac{E|X_1| I\{|X_1|\le
H^{-1}(d/n)\}}{(\alpha/(1-\alpha)) H^{-1}(d/n) d/n}=1.
\]
Similarly,
\begin{eqnarray*}
&&E X_1^2 I\{|X_1| \le
H^{-1}(d/n)\}\\&&\quad=\int_{-H^{-1}(d/n)}^{H^{-1}(d/n)} x^2\, \mathrm{d}F(x)\\
&&\quad=-\int_0^{H^{-1}(d/n)} x^2\,\mathrm{d}H(x)=-x^2 H(x) |_{H^{-1}(d/n)}+2
\int_0^{H^{-1}(d/n)}x H(x)\,\mathrm{d}x.
\end{eqnarray*}
Again using \cite{binghametal}, page~26, we conclude
that
\[
\lim_{n\to\infty} \frac{E X_1^2 I\{|X_1|\le
H^{-1}(d/n)\}}{(H^{-1}(d/n))^2 d/n}=\frac{\alpha}{2-\alpha}.
\]
Hence, \eqref{eq:V_var} is established.

Arguing as above, we get
\begin{eqnarray*}
EX_1^4 I\{|X_1|\le H^{-1}(d/n)\}
&=&-\int_0^{H^{-1}(d/n)} x^4 \,\mathrm{d}H(x)
\\
&=&-x^4 H(x)|_{H^{-1}(d/n)} + 4 \int_0^{H^{-1}(d/n)} x^3 H(x)\, \mathrm{d}x
\end{eqnarray*}
and therefore
\[
\lim_{n\to\infty} \frac{E X_1^4 I\{|X_1|\le
H^{-1}(d/n)\}}{(H^{-1}(d/n))^4 d/n} = \frac{\alpha}{4-\alpha}.
\]
Similarly,
\[
\lim_{n\to\infty} \frac{E |X_1|^3 I\{|X_1|\le
H^{-1}(d/n)\}}{(H^{-1}(d/n))^3 d/n} =\frac{\alpha}{3-\alpha},
\]
completing the proof of \eqref{eq:4_mom}.
\end{pf}

\begin{pf*}{Proof of Theorem~\ref{th:weak_conv_V}}
Clearly, for each $n$, $X_j I\{|X_j|\le H^{-1}(d/n)\}$, $1\le j \le n$,
are independent and identically distributed random variables. By
Lemma~\ref{le:2.1}, we have
\[
\lim_{n\to\infty} \frac{ \sum_{j=1}^n E[X_j I\{|X_j|\le
H^{-1}(d/n)\}-E [X_j I\{|X_j|\le H^{-1}(d/n)\}]
]^4}{\bigl(\sum_{j=1}^n \var(X_j I\{|X_j| \le H^{-1}(d/n)\})\bigr)^2}=0,
\]
so the Lyapunov condition is satisfied. Hence, the result follows
immediately from Skorokhod~\cite{sko}.
\end{pf*}

A series of lemmas is needed to establish
Theorem~\ref{th:diff_trim_trunc}. Let $\eta_{n,1}\ge \eta_{n,2}\ge
\cdots\ge \eta_{n,n}$ denote the order statistics of
$|X_1|,\ldots,|X_n|$, starting with the largest value.

\begin{Lemma}\label{le:2.2}
If $H_0$ and \eqref{eq:F_cont} hold, then
\[
\{H(\eta_{n,k}),1\le k\le n\}\stackrel{\mathcal{D}}{=}\{S_k/S_{n+1},1\le
k\le n\},
\]
where
\[
S_k=e_1+\cdots+e_k, \qquad 1\le k\le n,
\]
and $e_1,e_2,\ldots,e_{n+1}$ are independent, identically distributed
exponential random variables with $E e_j=1$.
\end{Lemma}

\begin{pf}
The representation in Lemma~\ref{le:2.2} is well known (see, e.g.,
Shorack and Wellner \cite{showe}, page 335).
\end{pf}

Let $\eta_{n,d}(j)$ denote the $d$th largest among
$|X_1|,\ldots,|X_{j-1}|,|X_{j+1}|,\ldots,|X_n|$.

\begin{Lemma}\label{le:2.3}
If $H_0$, \eqref{eq:ass_stable}, \eqref{eq:ass_symm}, \eqref{eq:dn},
\eqref{eq:F_cont} and \eqref{eq:mn} hold, then
\[
\sum_{j=1}^n \bigl|X_j \bigl(I\{|X_j|\le \eta_{n,d}\}-I\{|X_j|\le
\eta_{n,d}(j)\}\bigr)\bigr|=\mathrm{o}_P(A_n).
\]
\end{Lemma}

\begin{pf}
First, we note that $\eta_{n,d}(j)=\eta_{n,d}$ or
$\eta_{n,d}(j)=\eta_{n,d+1}$. Hence,
\[
\frac{H(\eta_{n,d})}{H(\eta_{n,d}(j))}\ge
\frac{H(\eta_{n,d})}{H(\eta_{n,d+1})}.
\]
By Lemma~\ref{le:2.2} and the law of large numbers, we have
\[
\frac{H(\eta_{n,d})}{H(\eta_{n,d+1})}\stackrel{\mathcal{D}}{=}
\frac{S_{d}}{S_{d+1}}=\frac{S_d}{S_d+e_{d+1}}=\frac{1}{1+e_{d+1}/S_d}=1+\mathrm{O}_P(d^{-1}).
\]
Furthermore, by the central limit theorem, we conclude that
\[
S_r=r\bigl(1+\mathrm{O}_P(r^{-1/2})\bigr)
\]
and thus
\[
H(\eta_{n,d})=\frac{d}{n}\bigl(1+\mathrm{O}_P(d^{-1/2})\bigr).
\]
Hence, for every $\varepsilon>0$, there is a constant
$C=C(\varepsilon)$ and an event $A=A(\varepsilon)$ such that $P(A)\ge
1-\varepsilon$, and on $A$,
\begin{equation}\label{eq:H_eta_qu}
\frac{H(\eta_{n,d})}{H(\eta_{n,d+1})}\ge 1-\frac{C}{d}
\end{equation}
and
\begin{equation}\label{eq:H_eta}
H(\eta_{n,d})\ge \frac{d}{n} \biggl(1-\frac{C}{\sqrt{d}}\biggr).
\end{equation}
We note that  $H(|X_j|)$ is uniformly distributed on $[0,1]$ and is
independent of $\eta_{n,d}(j)$. So, using \eqref{eq:H_eta_qu} and
\eqref{eq:H_eta}, we obtain that
\begin{eqnarray*}
&&
E \bigl[\bigl|X_j \bigl(I\{|X_j|\le \eta_{n,d}\}-I\{|X_j|\le\eta_{n,d}(j)\}\bigr)\bigr| I\{A\} \bigl]
\\
&&\quad=
E[|X_j|I\{\eta_{n,d}(j)\le |X_j| \le \eta_{n,d}\} I\{A\}]
\\
&&\quad\le
H^{-1}\biggl(\frac{d}{n}\biggl(1-\frac{C}{\sqrt{d}}\biggr)\biggr) E[I\{H(\eta_{n,d})\le H(|X_j|)\le H(\eta_{n,d}(j))\} I\{A\}]
\\
&&\quad\le
H^{-1}\biggl(\frac{d}{n}\biggl(1-\frac{C}{\sqrt{d}}\biggr)\biggr)EI\biggl\{H(\eta_{n,d}(j))\biggl(1-\frac{C}{d}\biggr)\le H(|X_j|)\le H(\eta_{n,d}(j))\biggr\}
\\
&&\quad\le
H^{-1}\biggl(\frac{d}{n}\biggl(1-\frac{C}{\sqrt{d}}\biggr)\biggr) E H(\eta_{n,d}(j))\frac{C}{d} \le H^{-1}\biggl(\frac{d}{n}\biggr(1-\frac{C}{\sqrt{d}}\biggr)\biggr)\frac{d+1}{n+1} \frac{C}{d}
\end{eqnarray*}
since $H(\eta_{n,d}(j))\le H(\eta_{n,d+1})$ and, by Lemma~\ref{le:2.2},
we have $EH(\eta_{n,d+1})=(d+1)/(n+1)$. The slow variation and
monotonicity of $H^{-1}$ yield
\[
\lim_{n\to\infty} \frac{
H^{-1}((d/n)(1-C/\sqrt{d}))}{ H^{-1}(d/n)}=1
\]
and thus we get that
\[
\lim_{n\to\infty} \frac{1}{A_n}\sum_{j=1}^n E\bigl|X_j\bigl(I\{|X_j|\le
\eta_{n,d}\}-I\{|X_j|\le \eta_{n,d}(j)\}\bigr)\bigr| I\{A\}=0.
\]
Since we can choose $\varepsilon>0$ as small as we wish,
Lemma~\ref{le:2.3} is proved.
\end{pf}

\begin{Lemma}\label{le:2.4}
If the conditions of Lemma~\ref{le:2.3} are satisfied, then
\[
\frac{1}{A_n} \sum_{j=1}^n |m(\eta_{n,d})-m(\eta_{n,d}(j))|=\mathrm{o}_P(1).
\]
\end{Lemma}

\begin{pf}
This can be proven along the lines of the proof of Lemma~\ref{le:2.3}.
\end{pf}

Let
\[
\xi_j=X_j\bigl(I\{|X_j|\le \eta_{n,d}(j)\}-I\{|X_j|\le
H^{-1}(n/d)\}\bigr)-m(\eta_{n,d}(j)).
\]

\begin{Lemma}\label{le:2.5}
If the conditions of Theorem~\ref{th:weak_conv_T} are satisfied, then
there is an $a>0$ such that for all $\tau>1/\alpha$ and
$0<\varepsilon<1/2$,
\begin{eqnarray}
\label{eq:E_xi}E \xi_j
&=&
0,
\\
\label{eq:E_xi2}E\xi_j^2
&=&
E\xi_1^2=\mathcal{O}\bigl(\bigl(H^{-1}(d/n)\bigr)^2(d^{1/2+\varepsilon}/n)+n^{2\tau} \exp(-a d^{2\varepsilon})\bigr),
\\
\label{eq:cov_xi}E\xi_i\xi_j
&=&
E\xi_1\xi_2 =\mathcal{O}\bigl(\bigl(H^{-1}(d/n)\bigr)^2(d^{1/2+3\varepsilon}/n^2)+n^{2\tau} \exp(-a d^{2\varepsilon})\bigr)
\end{eqnarray}
for $1\le j\le n$ and $1\le i<j\le n,$ respectively.
\end{Lemma}

\begin{pf}
It follows from the independence of $X_j$ and $\eta_{n,d}(j)$ that
\[
E\xi_j=E(E(\xi_j|\eta_{n,d}(j)))=E\bigl(m(\eta_{n,d}(j))-m(\eta_{n,d}(j))\bigr)=0,
\]
so \eqref{eq:E_xi} is proved.

The first relation in \eqref{eq:E_xi2} is clear. For the second part,
we note that
\[
E \xi_1^2\le 2EX_1^2 \bigl(I\{|X_1|\le \eta_{n,d}(1)\}-I\{|X_1|\le
H^{-1}(d/n)\}\bigr)^2 + 2 E m^2(\eta_{n,d}(1))
\]
and
\begin{eqnarray*}
&&
E X_1^2\bigl(I\{|X_1|\le \eta_{n,d}(1)\}-I\{|X_1|\le
H^{-1}(d/n)\}\bigr)^2\\
&&\quad\le E X_1^2 I\{\eta_{n,d}(1)\le |X_1| \le H^{-1}(d/n)\} + E X_1^2
I\{H^{-1}(d/n)\le |X_1| \le \eta_{n,d}(1)\}\\
&&\quad\le \bigl(H^{-1}(d/n)\bigr)^2 P\{\eta_{n,d}(1)\le |X_1| \le
H^{-1}(d/n)\} \\
&&\qquad{} + E \bigl((\eta_{n,d}(1))^2 I\{H(\eta_{n,d}(1))\le H(|X_1|)\le
d/n\}\bigr).
\end{eqnarray*}
There are constants $c_1$ and $c_2$ such that
\begin{equation}\label{eq:P_S_k}
P\bigl\{|S_d-d|\ge x \sqrt{d}\bigr\} \le \exp(-c_1 x^2) \qquad \mbox{if }0\le x\le c_2
d.
\end{equation}
Let $0<\varepsilon<1/2$. Using Lemma~\ref{le:2.2} and \eqref{eq:P_S_k},
there is a constant $c_3$ such that
\begin{equation}\label{eq:PA}
P(A)\ge 1-c_3 \exp(-c_1 d^{2\varepsilon}),
\end{equation}
where
\[
A=\biggl\{\omega\dvtx \frac{d}{n} \biggl(1-\frac{1}{d^{1/2-\varepsilon}}\biggr)\le
H(\eta_{n,d}(1))\le \frac{d}{n} \biggl(1+\frac{1}{d^{1/2-\varepsilon}}\biggr)\biggr\}.
\]
Let $A^c$ denote the complement of $A$. By \eqref{eq:PA}, we have
\begin{eqnarray*}
&&\bigl(H^{-1}(d/n)\bigr)^2 P\{\eta_{n,d}(1)\le |X_1| \le H^{-1}(d/n)\}
\\
&&\quad= \bigl(H^{-1}(d/n)\bigr)^2 \bigl( P(A^c)+P\{\eta_{n,d}(1)\le |X_1|\le H^{-1}(d/n),
A\}\bigr)
\\
&&\quad\le \bigl(H^{-1}(d/n)\bigr)^2 \biggl(c_3 \exp(-c_1 d^{2\varepsilon})+P\biggl\{\frac{d}{n}\le
H(|X_1|)\le\frac{d}{n}\biggl(1+\frac{1}{d^{1/2-\varepsilon}}\biggr)\biggr\}\biggr)\\
&&\quad=\mathcal{O}\biggl(\bigl(H^{-1}(d/n)\bigr)^2 \biggl(\exp(-c_1
d^{2\varepsilon})+\frac{d^{1/2+\varepsilon}}{n}\biggr) \biggr).
\end{eqnarray*}
Similarly, by the independence of $|X_1|$ and $\eta_{n,d}(1)$, we have
\begin{eqnarray*}
&&E\bigl((\eta_{n,d}(1))^2 I\{H(\eta_{n,d}(1))\le H(|X_1|)\le d/n\}\bigr)\\
&&\quad\le E(\eta_{n,d}^2 (1) I\{A^c\})
\\
&&\qquad{} + E\bigl(\bigl(H^{-1}\bigl(d/n(1-d^{\varepsilon-1/2})\bigr)\bigr)^2
I\{d/n(1-d^{\varepsilon-1/2}) \le H(|X_1|)\le
d/n\}\bigr)\\
&&\quad=E(\eta_{n,d}^2(1) I\{A^c\})+\bigl(H^{-1}\bigl(d/n(1-d^{\varepsilon-1/2})\bigr)\bigr)^2
\frac{d}{n} d^{\varepsilon-1/2}.
\end{eqnarray*}
Since $H^{-1}(t)$ is a regularly varying function at 0 with index
$-1/\alpha$, for any $\tau>1/\alpha$, there is a constant $c_4$ such
that
\begin{equation}\label{eq:H_inv_upper}
H^{-1}(t)\le c_4 t^{-\tau},\qquad 0< t\le 1.
\end{equation}
By the Cauchy--Schwarz inequality, we have
\[
E\eta_{n,d}^2(1) I\{A^c\} \le (E \eta_{n,d}^4(1))^{1/2} (P(A^c))^{1/2}
\le (E \eta_{n,d}^4(1))^{1/2} c_3^{1/2} \exp\biggl(-\frac{c_1}{2}
d^{2\varepsilon}\biggr).
\]
Next, we use \eqref{eq:H_inv_upper} and Lemma~\ref{le:2.2} to conclude
that
\begin{eqnarray}\label{eq:eta_4_mom}
E \eta_{n,d}^4(1)&\le& E\eta_{n,d}^4\le c_4^4
E\biggl(\frac{S_d}{S_{n+1}}\biggr)^{-4\tau}
=c_4^4 E\biggl(1+\frac{S_{n+1}-S_d}{S_d}\biggr)^{4\tau} \nonumber\\[-8pt]\\[-8pt]
&\le& c_5 \biggl(1+E(S_{n+1}-S_d)^{4\tau} E\frac{1}{S_{d}^{4\tau}}\biggr) \le c_6
n^{4\tau} \nonumber
\end{eqnarray}
since $S_{d}$ has a Gamma distribution with parameter $d$ and therefore
$ES_{d}^{-4\tau}<\infty$ if $d\ge d_0(\tau)$. Thus, we have that
\begin{eqnarray*}
&&E X_1^2\bigl(I\{|X_1|\le \eta_{n,d}(1)\}-I\{|X_1|\le H^{-1}(d/n)\}\bigr)^2\\
&&\quad= \mathcal{O}\biggl(\bigl(H^{-1}(d/n)\bigr)^2 (d^{\varepsilon+1/2}/n)+n^{2\tau}
\exp\biggl(-\frac{c_1}{2} d^{2\varepsilon}\biggr)\biggr).
\end{eqnarray*}
Similar arguments give
\[
E m^2(\eta_{n,d}(1))=\mathcal{O}\biggl(\bigl(H^{-1}(d/n)\bigr)^2
(d^{\varepsilon+1/2}/n)+n^{2\tau} \exp\biggl(-\frac{c_1}{2}
d^{2\varepsilon}\biggr)\biggr).
\]
The proof of \eqref{eq:E_xi2} is now complete.

The first relation of \eqref{eq:cov_xi} is trivial. To prove the second
part, we introduce $\eta_{n,d}(1,2)$, the $d$th largest among
$|X_3|,|X_4|,\ldots,|X_n|$. Set
\[
\xi_{1,2}=X_1\bigl(I\{|X_1|\le \eta_{n,d}(1,2)\}-I\{|X_1|\le
H^{-1}(d/n)\}\bigr)-m(\eta_{n,d}(1,2))
\]
and
\[
\xi_{2,1}=X_2\bigl(I\{|X_2|\le \eta_{n,d}(1,2)\}-I\{|X_2|\le
H^{-1}(d/n)\}\bigr)-m(\eta_{n,d}(1,2)).
\]
Using the independence of $|X_1|$, $|X_2|$ and $\eta_{n,d}(1,2)$, we
get
\begin{equation}\label{eq:xi_cov_0}
E\xi_{1,2} \xi_{2,1}=0.
\end{equation}
Next, we observe that
\begin{eqnarray*}
\xi_1 \xi_2&=&X_1\bigl(I\{|X_1|\le \eta_{n,d}(1)\}-I\{|X_1|\le \eta_{n,d}(1,2)\} \xi_2\bigr)
\\
&&{}-\bigl(m(\eta_{n,d}(1))-m(\eta_{n,d}(1,2))\bigr) \xi_2\\
&&{} + X_2\bigl(I\{|X_2|\le \eta_{n,d}(2)\}-I\{|X_2|\le \eta_{n,d}(1,2)\}\bigr)
\xi_{1,2}
\\
&&{}-\bigl(m(\eta_{n,d}(2))-m(\eta_{n,d}(1,2))\bigr)
\xi_{1,2} +\xi_{1,2} \xi_{2,1}.
\end{eqnarray*}
So, by \eqref{eq:xi_cov_0}, we have
\begin{eqnarray*}
E \xi_1 \xi_2 &=& E\bigl(X_1 I\{\eta_{n,d}(1,2)<|X_1|\le \eta_{n,d}(1)\}
\xi_2\bigr)
+E\bigl(\bigl(m(\eta_{n,d}(1,2))-m(\eta_{n,d}(1))\bigr)\xi_2\bigr)\\
&&{}+ E\bigl(X_2 I\{\eta_{n,d}(1,2)< |X_2|\le \eta_{n,d}(2)\} \xi_{1,2}\bigr)
+E\bigl(\bigl(m(\eta_{n,d}(1,2))-m(\eta_{n,d}(2))\bigr)\xi_{1,2}\bigr)\\
&=&a_{n,1}+\cdots+a_{n,4}.
\end{eqnarray*}
It is easy to see that
\[
\eta_{n,d+2}\le \eta_{n,d}(1,2)\le \eta_{n,d}(1)\le \eta_{n,d}
\]
and
\[
\eta_{n,d+2}\le \eta_{n,d}(1,2)\le \eta_{n,d}(2)\le
\eta_{n,d}.
\]
Hence,
\[
\frac{H(\eta_{n,d}(1))}{H(\eta_{n,d}(1,2))}\geq\frac{H(\eta_{n,d})}{H(\eta_{n,d+2})}\stackrel{{\mathcal
D}}{=}\frac{S_d}{S_{d+2}}=1-\frac{e_{d+1}+e_{d+2}}{S_{d+2}},
\]
according to Lemma \ref{le:2.2}. Using (\ref{eq:P_S_k}), we get, for
any $0<\varepsilon<1/2,$
\[
P\bigl\{|S_{d+2}-(d+2)|\geq d^{2\varepsilon}\sqrt{d+2}\bigr\}\leq
\exp(-c_1d^{2\varepsilon}).
\]
The random variables $e_{d+1}$ and $e_{d+2}$ are exponentially
distributed with parameter 1 and therefore
\[
P\{e_{d+1}\geq d^{2\varepsilon}\}=P\{e_{d+2}\geq d^{2\varepsilon}\}\leq
\exp (-d^{2\varepsilon}).
\]
Thus,  for any $0<\varepsilon<1/2$, we obtain
\[
P\biggl\{\frac{H(\eta_{n,d}(1))}{H(\eta_{n,d}(1,2))}\ge 1-\frac{c_7
d^{2\varepsilon}}{d}\biggr\}\ge 1-c_8 \exp(-c_9 d^{2\varepsilon})
\]
and similar arguments yield
\[
P\biggl\{\frac{H(\eta_{n,d}(2))}{H(\eta_{n,d}(1,2))}\ge 1-\frac{c_7
d^{2\varepsilon}}{d}\biggr\}\ge 1-c_8 \exp(-c_9 d^{2\varepsilon})
\]
and
\[
P\biggl\{\frac{d}{n} \biggl(1-\frac{1}{d^{1/2-\varepsilon}}\biggr) \le H(\eta_{n,d})\le
\frac{d}{n} \biggl(1+\frac{1}{d^{1/2-\varepsilon}}\biggr)\biggr\} \ge 1-c_8 \exp(-c_9
d^{2\varepsilon})
\]
with some constants $c_7$, $c_8$ and $c_9$. We now define the event $A$
as the set on which
\[
\frac{H(\eta_{n,d}(1))}{H(\eta_{n,d}(1,2))}\ge
1-\frac{c_7}{d^{1-2\varepsilon}},\qquad
\frac{H(\eta_{n,d}(2))}{H(\eta_{n,d}(1,2))}\ge
1-\frac{c_7}{d^{1-2\varepsilon}}
\]
and
\[
\frac{d}{n} \biggl(1-\frac{1}{d^{1/2-\varepsilon}}\biggr) \le H(\eta_{n,d})\le
\frac{d}{n} \biggl(1+\frac{1}{d^{1/2-\varepsilon}}\biggr)
\]
hold. Clearly,
\[
P(A^c)\le 3 c_8 \exp(-c_9 d^{2\varepsilon}).
\]
Using the definition of $\xi_2$, we get that
\begin{eqnarray*}
a_{n,1}&\le& E\bigl(|X_1| I\{\eta_{n,d}(1,2)\le |X_1| \le
\eta_{n,d}(1)\}
\\
&&\quad{}\times |X_2|\bigl|I\{|X_2|\le \eta_{n,d}(2)\}-I\{|X_2|\le
H^{-1}(n/d)\}\bigr|\bigr)
\\
&&{}+E|X_1| I\{\eta_{n,d}(1,2)\le |X_1| \le \eta_{n,d}(1)\}
|m(\eta_{n,d}(2))|
\\
&\le& E |X_1| |X_2| I\{\eta_{n,d}(1,2)\le |X_1|\le \eta_{n,d}(1)\}
I\{H^{-1}(d/n)\le |X_2|\le \eta_{n,d}(2)\}
\\
&&{}+E|X_1| |X_2| I\{\eta_{n,d}(1,2)\le |X_1| \le
\eta_{n,d}(1)\} I\{\eta_{n,d}(2)\le |X_2|\le H^{-1}(d/n)\}
\\
&&{}+E|X_1| I\{\eta_{n,d}(1,2)\le |X_1|\le \eta_{n,d}(1)\}
|m(\eta_{n,d}(2))|
\\
&=&a_{n,1,1}+a_{n,1,2}+a_{n,1,3}.
\end{eqnarray*}
Using the definition of $A$, we obtain that
\begin{eqnarray*}
a_{n,1,1}&\le& E|X_1 X_2| I\{\eta_{n,d}(1,2)\le |X_1|\le\eta_{n,d}(1)\}
I\{H^{-1}(d/n)\le |X_2|\le \eta_{n,d}(2)\}I\{A\}
\\
&&{}+E|X_1 X_2| I\{\eta_{n,d}(1,2)\le |X_1|\le\eta_{n,d}(1)\}
I\{H^{-1}(d/n)\le |X_2|\le \eta_{n,d}(2)\}I\{A^c\}
\\
&\le& E\biggl(|X_1 X_2| I\biggl\{H(\eta_{n,d}(1,2))
\biggl(1-\frac{c_7}{d^{1-2\varepsilon}}\biggr) \le H(|X_1|)\le H(\eta_{n,d}(1,2))
\biggr\}
\\
&&\hphantom{E\biggl(}{}\times I\{A\} I\{H^{-1}(d/n)\le |X_2|\le \eta_{n,d}(2)\}\biggr)
+E(\eta_{n,d}^2 I\{A^c\})
\\
&\le&\biggl(H^{-1}\biggl(\frac{d}{n}\biggl(1-\frac{c_{10}}{d^{1/2-\varepsilon}}\biggr)\biggr)\biggr)^2
\\
&&\hphantom{\biggl(}{} \times E \biggl(I\biggl\{H(\eta_{n,d}(1,2))
\biggl(1-\frac{c_7}{d^{1-2\varepsilon}}\biggr) \le H(|X_1|) \le H(\eta_{n,d}(1,2))
\biggr\}
\\
&&\hspace*{10pt}\qquad{} \times I\biggl\{\frac{d}{n}\biggl(1-\frac{1}{d^{1/2-\varepsilon}}\biggr) \le
H(|X_2|) \le \frac{d}{n}\biggr\}\biggr) + E(\eta_{n,d}^2 I\{A^c\}).
\end{eqnarray*}
Again using the independence of $|X_1|$, $|X_2|$ and
$\eta_{n,d}(1,2)$, we conclude that
\begin{eqnarray*}
&&E\biggl(I\biggl\{H(\eta_{n,d}(1,2)) \biggl(1-\frac{c_7}{d^{1-2\varepsilon}}\biggr)\le
H(|X_1|)\le
H(\eta_{n,d}(1,2)) \biggr\}\\
&&\hspace*{-6pt}\qquad{}\times I\biggl\{\frac{d}{n} \biggl(1-\frac{1}{d^{1/2-\varepsilon}}\biggr) \le
H(|X_2|) \le
\frac{d}{n}\biggr\} \biggr)\\
&&\quad= E H( \eta_{n,d}(1,2)) \frac{c_7}{d^{1-2\varepsilon}} \frac{d}{n}
\frac{1}{d^{1/2-\varepsilon}}\le \frac{d}{n-1} \frac{c_7}{n}
\frac{1}{d^{1/2-3\varepsilon}}.
\end{eqnarray*}
The Cauchy--Schwarz inequality yields
\[
E(\eta_{n,d}^2 I\{A^c\})\le (E\eta_{n,d}^4)^{1/2} (P(A^c))^{1/2}
=\mathcal{O}\biggl(n^{2\tau} \exp\biggl(-\frac{c_9}{2} d^{2\varepsilon}\biggr)\biggr)
\]
for all $\tau>1/\alpha$ on account of \eqref{eq:eta_4_mom}. We thus
conclude
\[
a_{n,1,1}=\mathcal{O}\biggl( \bigl(H^{-1}(d/n)\bigr)^2
(d^{1/2+3\varepsilon}/n^2)+n^{2\tau}
\exp\biggl(-\frac{c_9}{2}d^{2\varepsilon}\biggr)\biggr).
\]
Similar, but somewhat simpler, arguments imply that
\[
a_{n,1,2}+a_{n,1,3}=\mathcal{O}\biggl(\bigl(H^{-1}(d/n)\bigr)^2
(d^{1/2+3\varepsilon}/n^2)+n^{2\tau}
\exp\biggl(-\frac{c_9}{2}d^{2\varepsilon}\biggr)\biggr),
\]
resulting in
\begin{equation}\label{eq:a_n_1}
a_{n,1}=\mathcal{O}\biggl(\bigl(H^{-1}(d/n)\bigr)^2
(d^{1/2+3\varepsilon}/n^2)+n^{2\tau}
\exp\biggl(-\frac{c_9}{2}d^{2\varepsilon}\biggr)\biggr).
\end{equation}
Following the lines of the proof of \eqref{eq:a_n_1}, the same rates
can be obtained for $a_{n,2}$ and $a_{n,3}$.
\end{pf}

\begin{Lemma}\label{le:2.6}
If the conditions of Theorem~\ref{th:weak_conv_T} are satisfied, then
\[
\frac{1}{A_n} \max_{1\le k\le n} \Biggl|\sum_{j=1}^k \xi_j\Biggr| =\mathrm{o}_P(1).
\]
\end{Lemma}

\begin{pf}
It is easy to see that for any $1\le \ell_1\le \ell_2\le n$, we have
\begin{eqnarray*}
E\Biggl(\sum_{j=\ell_1}^{\ell_2} \xi_j\Biggr)^2&=&(\ell_2-\ell_1+1)E\xi_1^2 +
(\ell_2-\ell_1)(\ell_2-\ell_1+1)E\xi_1 \xi_2
\\
&\leq& (\ell_2-\ell_1+1)(E\xi^2_1+nE\xi_1\xi_2).
\end{eqnarray*}
Lemma \ref{le:2.5} and (\ref{pityu}) yield
\[
E\xi^2_1\leq c_1\frac{A_n^2}{n}[d^{-1/2+\varepsilon}+n^{2\tau+1}\exp
(-ad^{2\varepsilon})]
\]
and
\[
E\xi_1\xi_2\leq
c_2\frac{A_n^2}{n^2}[d^{-1/2+3\varepsilon}+n^{2\tau+2}\exp
(-ad^{2\varepsilon})]
\]
for all $0<\varepsilon<1/6.$ Hence, we conclude that
\[
E\Biggl(\sum_{j=\ell_1}^{\ell_2} \xi_j\Biggr)^2\leq
c_3(\ell_2-\ell_1+1)\frac{A_n^2}{n}[d^{-1/2+3\varepsilon}+n^{2\tau+2}\exp
(-ad^{2\varepsilon})].
\]
So, using an inequality of Menshov (see Billingsley \cite{billingsley}, page 102),
we get that
\begin{eqnarray*}
E\Biggl(\max_{1\le k\le n} \Biggl|\sum_{j=1}^k \xi_j\Biggr|\Biggr)^2&\leq& c_4(\log n)^2
A_n^2[d^{-1/2+3\varepsilon}+n^{2\tau+2}\exp (-ad^{2\varepsilon})]
\\
&\leq& c_4A_n^2\bigl[(\log n)^2d^{-2/7}+\exp \bigl((2\tau+2)\log n+ 2\log\log
n-ad^{2\varepsilon}\bigr)\bigr]
\\
&=&A_n^2\mathrm{o}(1)\qquad\mbox{as } n\rightarrow \infty,
\end{eqnarray*}
where $\varepsilon=1/14$ and $d=(\log n)^\gamma$ with any $\gamma
>7$, resulting in
\[
\frac{1}{A_n^2} E\Biggl(\max_{1\le k\le n} \Biggl|\sum_{j=1}^k \xi_j\Biggr|\Biggr)^2=\mathrm{o}(1).
\]
Markov's inequality now completes  the proof of Lemma~\ref{le:2.6}.
\end{pf}

\begin{pf*}{Proof of Theorem~\ref{th:diff_trim_trunc}}
This follows immediately from Lemmas~\ref{le:2.3}, \ref{le:2.4} and
\ref{le:2.6}.
\end{pf*}

\begin{pf*}{Proof of Theorem~\ref{th:weak_conv_T}}
According to \eqref{eq:star}, Theorems~\ref{th:weak_conv_V} and
\ref{th:diff_trim_trunc} imply Theorem~\ref{th:weak_conv_T}.
\end{pf*}

\begin{Lemma}\label{le:2.7}
If the conditions of Theorem~\ref{th:weak_conv_T} are satisfied, then
\[
\frac{\hat{A}_n}{A_n}\longrightarrow 1  \qquad\mbox{a.s.}
\]
\end{Lemma}

\begin{pf}
This is an immediate consequence of Haeusler and Mason \cite{hama}.
\end{pf}

\begin{pf*}{Proof of Theorem~\ref{th:weak_conv_T2}}
From Slutsky's lemma, it follows that Lemma~\ref{le:2.7} and
Theorem~\ref{th:weak_conv_T} imply the result.
\end{pf*}

\begin{pf*}{Proof of Example~\ref{ex:diff_trim_trunc}}
Since $H'(x)=-f(x)$,  our assumptions imply that  $H'(x)$  is also
regularly varying at $\infty$.  By elementary results on regular
variation (see, e.g., \cite{binghametal}), it follows
that
\[
H(x)=1-F(x)=\int_x^\infty f(t)\,\mathrm{d}t \sim \frac{1}{\alpha} xf(x) \qquad\mbox{as }
x\to\infty.
\]
Hence, $H^{-1}$ is regularly varying at $0$ and therefore the function
$(H^{-1}(t))'=1/H'(H^{-1}(t))$ is also regularly varying at $0$. Also,
\[
m'(x)=\frac{\mathrm{d}}{\mathrm{d}x}\int_0^x tf(t)\,\mathrm{d}t=xf(x)\sim {\alpha}H(x) \qquad\mbox{as }
x\to \infty
\]
and therefore $m'(H^{-1}(t))\sim t\alpha$. Using Lemma~\ref{le:2.2},
the mean value theorem gives
\[
\frac{nm(\eta_{n,d})}{B_n}\stackrel{\mathcal D
}{=}\frac{nm(H^{-1}(S_d/S_{n+1}))}{B_n}
=\frac{n(\ell(S_d/S_{n+1})-\ell(d/n))}{B_n}
 =\frac{n}{B_n}
\ell'(\xi_n) \biggl(\frac{S_d}{S_{n+1}}-\frac{d}{n}\biggr),
\]
where $\xi_n$ is between $S_d/S_{n+1}$ and $d/n$, and
$\ell(t)=m(H^{-1}(t)).$ It follows from the central limit theorem for
central order statistics that
\begin{equation}\label{reis}
\frac{n}{d^{1/2}} \biggl(\frac{S_d}{S_{n+1}}-\frac{d}{n}\biggr) \Dconv N(0,1).
\end{equation}
The regular variation of $\ell'$ and (\ref{reis}) yield
\[
\ell'(\xi_n) /\ell'(d/n) \to 1\qquad\mbox{in probability}.
\]
The result now follows from (\ref{reis}) by observing that
\[
\frac{n}{B_n}\ell'(d/n)\sim \frac{n}{d^{1/2}}.
\]
\upqed\end{pf*}

The proof of Theorem~\ref{th:weak_conv_T3} is based on analogs of
Theorems~\ref{th:weak_conv_V}, \ref{th:diff_trim_trunc} and
Lemmas~\ref{le:2.3}--\ref{le:2.7} when $E X_j^2<\infty$.

\begin{Lemma}\label{le:2.8}
If the conditions of Theorem~\ref{th:weak_conv_T3} are satisfied, then
\[
\frac{1}{\sqrt{n}} \sum_{j=1}^{\nt} \bigl(X_j I\{|X_j|\le H^{-1}(d/n)\}-E
[X_1 I\{|X_1|\le H^{-1}(d/n)\}]\bigr)
\stackrel{\mathcal{D}[0,1]}{\longrightarrow} \sigma W(t),
\]
where $\sigma^2=\var X_1$.
\end{Lemma}

\begin{pf}
By $E X_1^2<\infty$, we have
\[
E\bigl[X_1 I\{|X_1|\le H^{-1}(d/n)\}-E[X_1 I\{|X_1|\le
H^{-1}(d/n)\}]-(X_1-E X_1)\bigl]^2 \longrightarrow 0
\]
as $n\to \infty$. So, using L\'{e}vy's inequality \cite{loeve}, page 248, we get
\begin{eqnarray*}
&&\frac{1}{\sqrt{n}} \max_{1\le k\le n} \Biggl| \sum_{j=1}^k \bigl(X_j
I\{|X_j|\le H^{-1}(d/n)\}\\
&&\hphantom{\frac{1}{\sqrt{n}} \max_{1\le k\le n} \Biggl| \sum_{j=1}^k \bigl(}{}- E[X_1 I\{|X_1|\le H^{-1}(d/n)\}]-(X_j-E
X_1)\bigr)\Biggr|=\mathrm{o}_P(1).
\end{eqnarray*}
Donsker's theorem  (see \cite{billingsley}, page 137) now implies the result.
\end{pf}

\begin{Lemma}\label{le:2.9}
If the conditions of Theorem~\ref{th:weak_conv_T3} are satisfied, then
\[
\frac{1}{\sqrt{n}} \sum_{j=1}^n \bigl| X_j \bigl(I\{|X_j|\le
\eta_{n,d}\}-I\{|X_j| \le \eta_{n,d}(j)\}\bigr)\bigr| = \mathrm{o}_P(1)
\]
and
\[
\frac{1}{\sqrt{n}} \sum_{j=1}^n
|m(\eta_{n,d})-m(\eta_{n,d}(j))|=\mathrm{o}_P(1).
\]
\end{Lemma}

\begin{pf}
We adapt the proof of Lemma~\ref{le:2.3}. We recall that $A$ is an
event satisfying \eqref{eq:H_eta_qu}, \eqref{eq:H_eta} and $P(A)\ge
1-\varepsilon$, where $\varepsilon>0$ is an arbitrary small positive
number. We also showed that
\begin{eqnarray*}
&&E\bigl(\bigl|X_j \bigl(I\{|X_j|\le \eta_{n,d}\}-I\{|X_j|\le \eta_{n,d}(j)\}\bigr)\bigr|
I\{A\}\bigr)
\\&&\quad \le H^{-1}\biggl(\frac{d}{n}\biggl(1-\frac{C}{\sqrt{d}}\biggr)\biggr) \frac{d+1}{n+1}
\frac{C}{d}
\end{eqnarray*}
for some constant $C$. Assumption $E X_1^2<\infty$ yields
\[
\limsup_{x\to 0} x^{1/2} H^{-1}(x)<\infty
\]
and therefore
\[
\lim_{n\to\infty} \sqrt{n} H^{-1}\biggl(\frac{d}{n} \biggl(1-\frac{C}{\sqrt{d}}\biggr)\biggr)
\frac{d+1}{n+1} \frac{C}{d}=0
\]
for all $C>0$. Thus,  for all $\varepsilon>0$, we have
\[
\lim_{n\to\infty} \frac{1}{\sqrt{n}} \sum_{j=1}^n E\bigl|X_j \bigl(I\{|X_j|\le
\eta_{n,d}\}-I\{|X_j|\le \eta_{n,d}(j)\}\bigr)\bigr| I\{A\}=0.
\]
Since we can choose $\varepsilon>0$ as small as we wish, the first
result is proved. The second part of the lemma can be established
similarly.
\end{pf}

\begin{Lemma}\label{le:2.10}
If the conditions of Theorem~\ref{th:weak_conv_T3} are satisfied, then
for all $0<\varepsilon<1/2,$
\begin{eqnarray*}
E \xi_j&=&0,\qquad  1\le j\le n,\\
E\xi_j^2&=&E\xi_1^2=\mathcal{O}\bigl(\bigl(H^{-1}(d/n)\bigr)^2
d^{1/2+\varepsilon}/n+n \exp(-a d^{2\varepsilon})\bigr),\qquad1\le
j\le n,\\
E\xi_i \xi_j&=&E\xi_1\xi_2=\mathcal{O}\bigl(\bigl(H^{-1}(d/n)\bigr)^2
d^{1/2+3\varepsilon}/n^2+n \exp(-a d^{2\varepsilon})\bigr),\qquad 1\le i\neq
j\le n.
\end{eqnarray*}
\end{Lemma}

\begin{pf}
The proof of Lemma~\ref{le:2.5} can be repeated, only
\eqref{eq:H_inv_upper} should be replaced with
\begin{equation}\label{eq:Hinv_fin_var}
H^{-1}(t)\le C t^{-1/2},\qquad 0<t\le 1.
\end{equation}
\end{pf}

\begin{Lemma}\label{le:2.11}
If the conditions of Theorem~\ref{th:weak_conv_T3} are satisfied, then
\[
\frac{1}{\sqrt{n}} \max_{1\le k\le n} \Biggl|\sum_{j=1}^k \xi_j\Biggr| = \mathrm{o}_P(1).
\]
\end{Lemma}

\begin{pf}
Following the proof of Lemma~\ref{le:2.6}, we get
\begin{equation}
E \Biggl(\max_{1\le k\le n} \Biggl|\sum_{j=1}^k \xi_j\Biggr|\Biggr)^2 \le c_1n(\log
n)^2[d^{-1/2+3\varepsilon} +n^3\exp (-ad^{2\varepsilon})]=\mathrm{o}(n)
\end{equation}
as $n\to\infty$. Markov's inequality completes the proof of
Lemma~\ref{le:2.11}.
\end{pf}

\begin{Lemma}\label{le:2.12}
If the conditions of Theorem~\ref{th:weak_conv_T3} are satisfied, then
\[
\frac{1}{\sqrt{n}} \max_{1\le k\le n} \Biggl|\sum_{j=1}^k \bigl[X_j\bigl(I\{|X_j|\le
\eta_{n,d}\}-I\{|X_j|\le H^{-1}(d/n)\}\bigr)-m(\eta_{n,d})\bigr]\Biggr|=\mathrm{o}_P(1).
\]
\end{Lemma}

\begin{pf}
It follows immediately from Lemmas~\ref{le:2.9} and \ref{le:2.11}.
\end{pf}

\begin{pf*}{Proof of Theorem~\ref{th:weak_conv_T3}}
By Lemmas~\ref{le:2.8} and \ref{le:2.12}, we have that
\[
\frac{T_n(t)}{\sigma \sqrt{n}}
\stackrel{\mathcal{D}[0,1]}{\longrightarrow} B(t).
\]
It is easy to see that
\[
\frac{\hat{A}_n^2}{n} \stackrel{P}{\longrightarrow} \sigma^2,
\]
which completes the proof of Theorem~\ref{th:weak_conv_T3}.
\end{pf*}

\begin{pf*}{Proof of Theorem~\ref{th:resampling}}
We show that
\begin{equation}\label{eq:x_i_x_i2}
\frac{\max_{1\le j\le n} |x_j|}{\sqrt{\sum_{j=1}^n
x_j^2}}\longrightarrow 0 \qquad \mbox{a.s.}
\end{equation}
By Lemma~\ref{le:2.7} it is enough to prove that
\[
\frac{\max_{1\le j\le n}|x_j|}{A_n}\longrightarrow 0 \qquad \mbox{a.s.}
\]
It follows from the definition of $x_j$ that
\[
\max_{1\le j\le n}|x_j|\le \eta_{d,n}+|\bar{X}_{n,d}|\le 2 \eta_{d,n}.
\]
Using Kiefer \cite{kiefer71} (see Shorack and Wellner \cite{showe}), we
get
\[
\frac{\eta_{d,n}}{A_n}\longrightarrow 0 \qquad \mbox{a.s.}
\]
Since \eqref{eq:x_i_x_i2} holds for almost all realizations of
$X_1,X_2,\ldots,$ Theorem \ref{th:resampling} is implied by Ros\'{e}n \cite{ro} when we sample without replacement and by
Prohorov \cite{prohorov56} when we sample with replacement (bootstrap).
\end{pf*}

\begin{pf*}{Proof of Theorem~\ref{th:resampling_H_A}}
This can be established  along the lines of the proof of Theorem
\ref{th:resampling}.
\end{pf*}

\section*{Acknowledgements}

This research was supported by FWF
Grant S9603-N23 and OTKA Grants K 67961 and K 81928 (Istv\'{a}n
Berkes), partially supported by NSF Grant DMS-00905400 (Lajos
Horv\'{a}th) and partially supported by FWF Grant S9603-N23 (Johannes
Schauer).

\printhistory

\end{document}